\documentclass[10pt]{article}
\usepackage[letterpaper]{geometry}
\usepackage{hicss}
\usepackage{times}
\usepackage[none]{hyphenat}
\usepackage{url}
\usepackage{latexsym}
\usepackage{indentfirst}
\usepackage{graphicx}
\graphicspath{{images/}}

\usepackage{amsmath,amssymb,amsfonts}
\usepackage{algorithmic}
\usepackage{booktabs}
\newtheorem{theorem}{Theorem}[section]

\newtheorem{definition}[theorem]{Definition}
\newtheorem{corollary}[theorem]{Corollary}
\newtheorem{proposition}[theorem]{Proposition}

\usepackage{subfigure}
\usepackage{etoolbox}
 \usepackage[ruled,linesnumbered]{algorithm2e} 

\newtheorem{assumption}{Assumption}
\newcommand{\tabincell}[2]{\begin{tabular}{@{}#1@{}}#2\end{tabular}}
\def\BibTeX{{\rm B\kern-.05em{\sc i\kern-.025em b}\kern-.08em
   T\kern-.1667em\lower.7ex\hbox{E}\kern-.125emX}}

\newenvironment{proof_def}
{\textit{Proof:} }
{$\square$}

\title{Mechanism Design for Efficient Nash Equilibrium in Oligopolistic Markets}

\author{Kaiying Lin \\
 Electrical Engineering \\
  Southeast University \\
 {\underline{linkaiying\_seu@qq.com}} \\\And
 Beibei Wang \\
 Electrical Engineering \\
 Southeast University \\
 {\underline{wangbeibei@seu.edu.cn} }\\\And 
 Pengcheng You \\
 Electrical and Computer Engineering\\
 Johns Hopkins University \\
 {\underline{pcyou@jhu.edu}} \\}

\begin{document}
\maketitle
\begin{abstract}

This paper investigates the efficiency loss in social cost caused by strategic bidding behavior of individual participants in a supply-demand balancing market, and proposes a mechanism to fully recover equilibrium social optimum via subsidization and taxation. We characterize the competition among supply-side firms to meet given inelastic demand, with linear supply function bidding and the proposed efficiency recovery mechanism. We show that the Nash equilibrium of such a game exists under mild conditions, and more importantly, it achieves the underlying efficient supply dispatch and the market clearing price that reflects the truthful system marginal production cost.
Further, the mechanism can be tuned to guarantee self-sufficiency, i.e., taxes collected counterbalance subsidies needed. 
Extensive numerical case studies are run to validate the equilibrium analysis, and we employ individual net profit and a modified version of Lerner index as two metrics to evaluate the impact of the mechanism on market outcomes by varying its tuning parameter and firm heterogeneity.


\end{abstract}

\section{Introduction}


Deregulated markets are designed to foster competition among participants by allowing them to bid freely, in the aim of driving efficient operation and investment. However, natural barriers to market entry and strategic dominance of existing participants often lead to oligopolistic markets that are not socially efficient in principle, e.g., electricity markets. There exists systematic loss of efficiency attributed to the strategic behavior of individual participants, which exploits the anticipation of their bidding impact to manipulate market clearing prices. Such ability to earn themselves extra benefits is referred to as \emph{market power} in economics.



The study of the efficiency loss due to the exercise of market power has been prevailing for decades, especially from a game theoretic perspective. 
It is claimed in \cite{2004Efficiency} that in a network resource allocation market with a proportional pricing mechanism, the aggregate surplus at the Nash equilibrium is no worse than a factor of $4\sqrt{2}-5$ times the optimal
aggregate surplus, thus bounding the efficiency loss by approximately $34\%$.
\cite{2011Parameterized} looks at a market with parameterized supply function bidding and uniform pricing, and provides a similar upper bound on efficiency loss at a symmetric equilibrium. 
To ease analysis, \cite{2015Demand} considers a simpler form of linear supply function bids and is able to characterize the general Nash equilibrium in such a game. Further, the maximal efficiency loss can also be bounded under mild conditions. 

There have been extensive efforts towards alleviating such loss of efficiency, from both academia and industry.
\cite{2010Autonomous} builds a quantity-bid based market/platform where the unique Nash equilibrium achieves social optimum. The billing mechanism is designed to align individual payoff with global utility. 
A similar idea is adopted in \cite{2012Advanced} that employs the Vickrey–Clarke–Groves (VCG) mechanism to bill participants. The rationale of the VCG mechanism is to make each participant responsible for the increment of the total system cost due to its participation, which also guarantees social optimum at the Nash equilibrium.
However, this class of mechanisms require sharing certain private information.
\cite{2012Vehicle} develops an interaction mechanism for aggregators to solicit regulation service from electric vehicles, which guarantees to achieve the global efficiency of a carefully designed problem at the Nash equilibrium. However, the game setup and the structure of strategy space are highly specific and difficult to generalize.
Additionally, the industry, electricity markets in particular, has also implemented policies to mitigate the impact of market power, such as market power test combined with price substitution \cite{reitzes2007review}, and posterior penalty or adjustment in the UK electricity market \cite{2004Efficiency} and Nord Pool \cite{daskalakis2009electricity}. Most of them are empirically designed and may need further theoretical support.

Our work contributes to this line of research by proposing a market mechanism that fully recovers social optimum at the Nash equilibrium among firms (suppliers) bidding to meet inelastic demand. The mechanism uses only the market quantities, e.g., bids, dispatch and prices, to compute and impose an additional monetary term, in the form of subsidy or tax, on each individual firm. Such a mechanism, when co-executed with market clearing rules, drives the Nash equilibrium to be efficient that achieves the minimum social cost.
Furthermore, the market can control a tuning parameter of the mechanism to maintain self-sufficiency, i.e., taxes collected counterbalance subsidies needed, thus redistributing surplus among firms. We use individual net profit and a modified version of Lerner index as metrics in case studies to evaluate the impact of the mechanism on market outcomes, in terms of varying tuning parameters and firm heterogeneity. 


The main contributions are three-fold:

\noindent
\emph{1}. We design an efficiency recovery mechanism, parameterized by a tuning constant, that provably drives the Nash equilibrium of a linear supply function bidding game among firms to be socially optimal. Indeed, the rationale underneath this mechanism design can be used in other market setups beyond linear supply function bidding.

\noindent
\emph{2}. We show that the tuning parameter can be set such that the market is self-sufficient in the presence of the efficiency recovery mechanism, and uncover the intrinsic tradeoff in selecting this parameter between market self-sufficiency and individual profitability.

\noindent
\emph{3}. We specifically develop two evaluation metrics of individual net profit and modified Lerner index for the mechanism and extensively study its impact on market outcomes with varying tuning parameters and firm heterogeneity.

The remainder of this paper is organized as follows. Section~2 introduces the preliminaries of the market setup and equilibria characterization. 
Section~3 proposes the efficiency recovery mechanism, while
Section~4 discusses its implementation and evaluation.
Section~5 presents case studies that validate our analysis and illustrate the impact of the mechanism.
In the end, Section~6 concludes.


\section{Preliminaries}

\subsection{System Model}

We consider a supply-demand balancing market where a set $\mathcal{N}:=\{1,2,\dots,|\mathcal{N}|\}$ of firms bid linear supply functions to meet given inelastic demand $d>0$.
A linear supply function takes the form of 
\begin{equation}
    q_i=b_ip 
\end{equation}
for each firm $i\in\mathcal{N}$, where $q_i$ denotes its supply and $p$ is the market clearing price. The linear coefficient $b_i\ge0$ therefore indicates the amount of supply incentivized by per unit price and is freely chosen by each individual firm. The larger $b_i$ is, the more a firm is willing to supply at price $p$ -- a natural interpretation as firms' price sensitivity. Without ambiguity, we may refer to $b_i$ as supply function bids from firms.
Suppose that firm $i$ incurs a production cost $C_i(q_i)$ for its supply $q_i$. We adopt general assumptions on the cost function that $C_i(\cdot):\mathbb{R}\rightarrow \mathbb{R}$ is continuously differentiable, convex, and strictly increasing with $C_i(0)=0$.

After collecting the supply function bids $B:=(b_i,i\in\mathcal{N})$ from all firms, the market matches the total supply with the demand $d$ 
\begin{equation}\label{eq:supply-demand_balance}
    \sum_{i} q_i = d
\end{equation}
by setting a clearing price at
\begin{equation}\label{eq:clearing_price}
    p(B)=\frac{d}{\sum_i{b_i}}
\end{equation}
To avoid degenerate situations where all firms make bids of 0 and supply-demand balance cannot be met, we make the following assumption:
\begin{assumption}
If all firms bid $b_i =0$, $\forall i\in\mathcal{N}$, these bids will be rejected by the market.
\end{assumption}
For notational convenience, we further define $B_{-i}:=(b_j,j\in\mathcal{N}\backslash\{i\})$ and $\sum_{-i}b_j:=\sum_{j\in\mathcal{N}\backslash\{i\}}b_j$.

\subsection{Competitive Equilibrium}

Consider a perfectly competitive market where all firms are rational price takers, i.e., they take market clearing prices as given and respond by making optimal bids to maximize their profit. In particular, given price $p$, each firm $i$ solves the following individual bidding problem for profit maximization:
\begin{eqnarray}\label{eq:price_taker_bidding}  
      \max_{b_{i}\ge0} && \pi_i(b_i;p):= pq_{i} - C_{i}(q_{i}) \\ \nonumber
      && \qquad \quad \ \ \ \,  = p^2 b_i - C_i(pb_i)  
\end{eqnarray}
where $\pi_i$ denotes the profit of firm $i$.

In such a market, a competitive equilibrium is defined as follows.
\begin{definition}
A competitive equilibrium among price-taking firms is a tuple $(B,p)$ that satisfies
\begin{enumerate}
    \item $b_i$ is optimal w.r.t. \eqref{eq:price_taker_bidding}, given $p$, for $\forall i\in\mathcal{N}$;
    \item supply matches demand, i.e., \eqref{eq:supply-demand_balance}.
\end{enumerate}
\end{definition}

It is shown in \cite{2015Demand} that through the analysis of the optimality KKT conditions of \eqref{eq:price_taker_bidding}, such a competitive equilibrium exists and can be characterized by
\begin{proposition}\label{prop:competitive_equilibrium}
The competitive equilibrium $(B,p)$ among price-taking firms exists. Moreover, it is efficient in the sense that it minimizes the social cost, defined as
\begin{subequations}\label{eq:CE_opt}
\begin{eqnarray}
\min_{q_i\ge0,i\in\mathcal{N}} && \sum_{i} C_i(q_i) \\
\mathrm{s.t.} && \sum_{i}q_{i}=d 
\end{eqnarray}
\end{subequations}
with $q_i = b_i p$, $i\in\mathcal{N}$.
\end{proposition}
If all the cost functions $C_i(\cdot)$ are strictly convex, there will exist a unique competitive equilibrium. 

\subsection{Nash Equilibrium}

We further consider an oligopolistic market where all firms are aware of the market clearing rule \eqref{eq:clearing_price} and bid strategically. In particular, each firm $i$ will anticipate the impact of its bidding decision on the clearing price and integrate such anticipation into the optimal bidding problem:
\begin{subequations}\label{eq:strategic_bidding}
\begin{eqnarray}\label{eq:strategic_bidding.a}
 \max_{b_{i}\geq 0} && \pi (b_{i};B_{-i})= pq_{i}-C_{i}(q_{i}) \\\nonumber 
&& \qquad \qquad \ \, \,  = p^2 b_i - C_i(pb_i) \\
\mathrm{s.t.} && \eqref{eq:clearing_price} \label{eq:strategic_bidding.b} 
\end{eqnarray}
\end{subequations}

A Nash equilibrium in such a market is where no firm has an incentive to deviate from the current status unilaterally, formally defined by 
\begin{definition}
A Nash equilibrium among strategic firms is a tuple $(B,p)$ that satisfies
\begin{enumerate}
    \item $b_i$ is optimal w.r.t. \eqref{eq:strategic_bidding}, given $B_{-i}$, for $\forall i\in\mathcal{N}$;
    \item supply matches demand, i.e., \eqref{eq:supply-demand_balance}.
\end{enumerate}
\end{definition}

\cite{2015Demand} has also characterized the existence and structure of such a Nash equilibrium, summarized below.

\begin{proposition}
If $|\mathcal{N}| \ge 3$ holds, the Nash equilibrium among strategic firms exists and is unique. Moreover, it solves the following convex optimization problem:
\begin{subequations}\label{eq:NE_opt}
\begin{eqnarray}
\min_{0\le q_i\le \frac{d}{2},i\in\mathcal{N}} && \sum_{i} D_i(q_i) \\
\mathrm{s.t.} && \sum_{i}q_{i}=d 
\end{eqnarray}
\end{subequations}
with
\begin{small}
\begin{equation}
    D_i(q_i) :=   \frac{d-q_i}{d-2q_i} C_i(q_i) - \int_0^{q_i} \frac{d}{(d-2x)^2} C_i(x)\mathrm{d}x \ .
\end{equation}
\end{small}
\end{proposition}
\eqref{eq:NE_opt} differs from \eqref{eq:CE_opt} mainly in the objective function where the constructed cost functions $D_i(\cdot)$ are strictly convex in $[0,\frac{d}{2})$, thus guaranteeing the uniqueness of the Nash equilibrium.


\section{Efficient Nash Equilibrium}

In this supply function bidding market, as we compare \eqref{eq:NE_opt} with \eqref{eq:CE_opt}, in general the Nash equilibrium likely deviates from the \emph{efficiency} in terms of the social cost, achieved by the competitive equilibrium. Such loss of efficiency from a social perspective is common in market competition due to the strategic behavior of participants. However, it degrades the role of a free market in fostering efficient trading and prevents the economy from prospering in the long run. We propose in this section a novel mechanism, co-executed with the supply function bidding mechanism, to recover efficiency of oligopolistic market operation, i.e., its Nash equilibrium, via extra revenue or expenditure imposed on each individual firm in the form of subsidy or tax.

\subsection{Efficiency Recovery Mechanism}

The efficiency recovery mechanism enforces an extra term $\Delta \pi_i$ of revenue, or expenditure if it is negative, on each firm $i$ such that its current objective function is explicitly
\begin{equation}
\begin{aligned}
        \pi_i^\dagger \ := \  & \pi + \Delta \pi_i \\
         \ = \ &  p q_i - C_i(q_i) + \Delta \pi_i \ ,
\end{aligned}
\end{equation}
where the market clearing rule remains the same for the original transaction part. We design $\Delta \pi_i$ in a way that interlaces with the supply function bidding and therefore affects the market competition.

In particular, we will exploit the existing market quantities to define the extra term as
\begin{equation}\label{eq:def_mechanism}
    \Delta \pi_i  =\frac{q_{i}^2}{2\sum_{-i}b_j}-\phi\frac{d^{2}}{2\sum_{-i}b_j} \ ,  
\end{equation}
where $\phi$ is a constant to be set by the market. Remarkably, the first term is associated with individual bids through their respective supply $q_i$, while the second term is independent of each firm $i$'s bidding decision $b_i$. 
In principle, the first term serves as the incentive to encourage supply when the willingness of other firms to supply is low, i.e., their aggregate bids $\sum_{-i}b_j$ are small. The second term is controlled by the market to achieve re-balance of surplus among firms via subsidy, $\Delta \pi_i >0$, or tax, $\Delta \pi_i <0$. Moreover, we specify the following assumption to avoid degenerate situations:
\begin{assumption}
In the case of $\sum_{-i} b_j = 0$ for any firm $i$, $\Delta \pi_i$ is set to a sufficiently large constant.
\end{assumption}

In the presence of the efficiency recovery mechanism, the optimal bidding problem of each strategic firm $i$
turns out
\begin{subequations}\label{eq:strategic_bidding_modified}
\begin{eqnarray}\label{eq:strategic_bidding_modified.a}
 \max_{b_{i}\geq 0} && \pi^\dagger (b_{i};B_{-i})\\\nonumber
 & = &  pq_{i}-C_{i}(q_{i}) +\Delta \pi_i \\ \nonumber
& = &  p^2 b_i - C_i(pb_i) + \frac{q_{i}^2}{2\sum_{-i}b_j}-\phi\frac{d^{2}}{2\sum_{-i}b_j} \\
\mathrm{s.t.} && \eqref{eq:clearing_price} \label{eq:strategic_bidding_modified.b} 
\end{eqnarray}
\end{subequations}

Our main result points to the characterization of a Nash equilibrium among strategic firms in such a case, as defined below.
\begin{definition}\label{def:NE_modified}
A Nash equilibrium among strategic firms in the presence of the efficiency recovery mechanism is a tuple $(B,p)$ that satisfies
\begin{enumerate}
    \item $b_i$ is optimal w.r.t. \eqref{eq:strategic_bidding_modified}, given $B_{-i}$, for $\forall i\in\mathcal{N}$;
    \item supply matches demand, i.e., \eqref{eq:supply-demand_balance}.
\end{enumerate}
\end{definition}
More precisely, we show the existence of such a Nash equilibrium and that it is indeed efficient in terms of the social cost, as summarized below.

\begin{theorem}\label{teo:main_NE}
If there are at least two firms with $q_i^*>0$ at the social optimum of \eqref{eq:CE_opt}, the Nash equilibrium among strategic firms in the presence of the efficiency recovery mechanism exists. Moreover, it is efficient that minimizes the social cost, i.e., it is an optimal solution to \eqref{eq:CE_opt}.
\end{theorem}

\begin{proof_def}
We first characterize the strategic bidding behavior of each individual firm $i$ out of the bidding problem \eqref{eq:strategic_bidding_modified}. By substituting \eqref{eq:clearing_price} into $\pi^\dagger$, we obtain its explicit expression
\begin{equation}
\begin{aligned}
        \pi^\dagger  =& \frac{d^2b_i}{(\sum_j b_j)^2}  - C_i\left(\frac{db_i}{\sum_j b_j}\right) + \frac{d^2b_i^2}{2(\sum_j b_j)^2\sum_{-i}b_j} \\ 
        &-\phi\frac{d^{2}}{2\sum_{-i}b_j} \ .
\end{aligned}
\end{equation}

In the case of $\sum_{-i} b_j = 0$, the profit of firm $i$ boils down to
\begin{equation}
    \pi^\dagger = \frac{d^2}{b_i} - C_i(d) + \mathrm{constant} \ ,
\end{equation}
which suggests that firm $i$ has an incentive to bid a positive $b_i$ yet as small as possible, in order to gain infinite profit. Note that the market will reject the bids if firm $i$ makes zero bid. Therefore, no equilibrium exists in this case.

In the case of $\sum_{-i} b_j >0$, its first-order derivative can be computed as
\begin{small}
\begin{equation}\label{eq:1st-order_derivative}
\begin{aligned}
    \frac{\mathrm{d}{\pi_{i}^{\dagger}}}{\mathrm{d}{b_{i}}}
     \ = \ & \frac{d^2(\sum_{-i} b_j - b_i)}{(\sum_j b_j)^3}-\frac{d\sum_{-i} b_j }{(\sum_j b_j)^2}C_{i}^{'}\left(\frac{db_i}{\sum_j b_j}\right) \\
     & + \frac{d^2b_i}{(\sum_j b_j)^3} \\
     \ = \  & \frac{d^2}{(\sum_j b_j)^2} \left( \frac{\sum_{-i} b_j}{\sum_j b_j}- \frac{\sum_{-i} b_j}{d} C_{i}^{'}\left(\frac{db_i}{\sum_j b_j}\right) \right)   \ .    
\end{aligned}
\end{equation}
\end{small}%
Inside the parenthesis, the first term, bounded by $(0,1]$, is decreasing in $b_i$, while the second term, bounded by $\left[\frac{\sum_{-i} b_j}{d} C_{i}^{'}\left(0\right), \frac{\sum_{-i} b_j}{d} C_{i}^{'}\left(d\right) \right)$, is non-decreasing in $b_i$. Note that the outer term is always positive. 
\begin{enumerate}
    \item[(a)] If $\frac{\sum_{-i} b_j}{d} C_{i}^{'}\left(0\right) \geq 1$ holds, we have $\frac{\mathrm{d}{\pi_{i}^{\dagger}}}{\mathrm{d}{b_{i}}} \leq 0$ for $\forall b_{i}\ge 0$ and $\frac{\mathrm{d}{\pi_{i}^{\dagger}}}{\mathrm{d}{b_{i}}} < 0$ for $\forall b_{i} > 0$. Therefore, $b_{i}=0$ maximizes the profit of firm $i$ and is the optimal bid.
    \item[(b)] If $\frac{\sum_{-i} b_j}{d} C_{i}^{'}\left(0\right) < 1$ holds, there exists a unique bid $0<b_i^*<\infty$ that satisfies 
    \begin{equation}\label{eq:unconstrained_opt}
         \frac{\sum_{-i} b_j}{\sum_j b_j} =  \frac{\sum_{-i} b_j}{d} C_{i}^{'}\left(\frac{db_i}{\sum_j b_j}\right)
    \end{equation}
    and thus $\frac{\mathrm{d}{\pi_{i}^{\dagger}}}{\mathrm{d}{b_{i}}} = 0$. For $0\le b_{i} < b_i^*$,  $\frac{\mathrm{d}{\pi_{i}^{\dagger}}}{\mathrm{d}{b_{i}}}>0$ holds, while for $b_i^* < b_{i}<\infty$, $\frac{\mathrm{d}{\pi_{i}^{\dagger}}}{\mathrm{d}{b_{i}}}<0$ holds. Therefore, $b_i^*$ maximizes the profit of firm $i$ and is the optimal bid.
\end{enumerate}

Note that the above two optimality conditions can be equivalently summarized as: for $\forall b_{i} \ge 0$,
\begin{equation}\label{eq:opt_summary}
   \left(\frac{d}{\sum_{-i}b_j + b_i^*}-C_{i}^{'}\left(\frac{db^*_i}{\sum_{-i}b_j + b_i^*}\right)\right) (b_{i}-b_{i}^{*})\leq 0 \ , 
\end{equation}
given $\sum_{-i}b_j$. If $b_i^* >0$ holds, \eqref{eq:opt_summary} enforces \eqref{eq:unconstrained_opt}, i.e., the situation (b). If $b_i^* =0$ holds, \eqref{eq:opt_summary} enforces $\frac{\sum_{-i} b_j}{d} C_{i}^{'}\left(0\right) \geq 1$, i.e., the situation (a).

By Definition~\ref{def:NE_modified}, the equilibrium conditions boil down to 
\begin{small}
\begin{subequations}\label{eq:eq_conditions}
\begin{eqnarray}\label{eq:eq_conditions.a}
   & \left(p-C_{i}^{'}\left(b^*_i p \right)\right) (b_{i}-b_{i}^{*})\leq 0  , \   \forall b_i \ge 0  , \ \forall i \in\mathcal{N} ,  \\\label{eq:eq_conditions.b}
   & \sum_i b^*_i p = d  \  ,
\end{eqnarray}
\end{subequations}
\end{small}%
with $p>0$.
Substituting the supply function form $q=b_ip$ and multiplying \eqref{eq:eq_conditions.a} with $p$, we equivalently arrive at
\begin{small}
\begin{subequations}\label{eq:social_opt_conditions}
\begin{eqnarray}\label{eq:social_opt_conditions.a}
   & \left(p-C_{i}^{'}\left(q_i^* \right)\right) (q_{i}-q_{i}^{*})\leq 0  , \   \forall q_i \ge 0  , \ \forall i \in\mathcal{N} ,  \\\label{eq:social_opt_conditions.b}
   & \sum_i q^*_i = d  \  ,
\end{eqnarray}
\end{subequations}
\end{small}%
which is exactly the KKT optimality conditions for the social cost minimization problem \eqref{eq:CE_opt}. Note that the prerequisite $\sum_{-i}b_j >0$ for $\forall i\in\mathcal{N}$ necessitates that there should be at least two firms with $q_i^*>0$ at the optimum of \eqref{eq:CE_opt}. Under such a circumstance, the Nash equilibrium solves \eqref{eq:CE_opt}, always exists, and is efficient. \end{proof_def}

In Theorem~\ref{teo:main_NE}, the condition that requires at least two active firms with $q_i^*>0$ at the social optimum suggests that in the presence of only one dominant firm in the market, it is difficult to restrict its strategic behavior which may arbitrarily manipulate the clearing price.


\subsection{Subsidy or Tax?}

Notably, this efficiency recovery mechanism in fact represents a family of functions $\Delta \pi_i$, parameterized by the constant $\phi$. The market can control $\phi$ to reallocate surplus among firms via subsidy, $\Delta \pi_i >0$, or tax, $\Delta \pi_i <0$. 
Given a $\phi$, whether a firm enjoys subsidy or incurs tax at the Nash equilibrium can be explicitly deduced as
\begin{corollary}
At the Nash equilibrium, for each strategic firm $i$, $\Delta \pi_i$ is a subsidy, if $\phi<\left(\frac{q_i^*}{d}\right)^2$ holds, while a tax if $\phi>\left(\frac{q_i^*}{d}\right)^2$ holds, where $(q_i^*,i\in\mathcal{N})$ is an optimal solution to \eqref{eq:CE_opt}.
\end{corollary}
The corollary is an immediate result of Theorem~\ref{teo:main_NE} and the definition of $\Delta \pi_i$ in \eqref{eq:def_mechanism}. Moreover, it indicates the following:
\begin{enumerate}
    \item[(a)] In the case of $\phi\le 0 \le \left(\frac{\inf_i q_{i}^*}{d}\right)^2$, no firm will incur tax at the Nash equilibrium;
     \item[(b)] In the case of $\phi\ge 1 \ge \left(\frac{\sup_i q_{i}^*}{d}\right)^2$, no firm will enjoy subsidy at the Nash equilibrium. 
\end{enumerate}
Therefore, a reasonable region for $\phi$ should be $[0,1)$: $\phi=0$ provides the minimum aggregate subsidy with every firm subsidized; $\phi=1$ instead overtaxes all firms that renders them unprofitable. Indeed, there exists a $\phi\in(0,1)$ to make the market self-sufficient, i.e., taxes collected counterbalance subsidies needed.
\begin{corollary}
By setting 
\begin{equation}\label{eq:self-sufficient_phi}
    \phi =  \frac{\sum_i^{}{\frac{{q_{i}^{*}}^2}{d-q_i^{*}}}}{d^2\sum_i^{}{\frac{1}{d - q_i^*}}} \  ,
\end{equation}
the market achieves self-sufficiency at the Nash equilibrium with 
\begin{equation}
    \sum_i{\Delta \pi_i }=0 \ .
\end{equation}
\end{corollary}
However, this particular $\phi$ in \eqref{eq:self-sufficient_phi} depends on the underlying social optimum $(q_i^*,i\in\mathcal{N})$ of \eqref{eq:CE_opt}, which is unknown \emph{a priori} since the private individual cost functions are not available. The setting of $\phi$ in real time, if it is allowed time-varying, can be potentially realized through online learning, and we would like to leave it as an open problem.

\section{Implementation and Evaluation}

In this section, we briefly discuss methods to reach the Nash equilibrium and evaluation metrics for the proposed efficiency recovery mechanism.

\subsection{Best Response Bidding Algorithm}

We employ a sequential best response bidding algorithm to reach the Nash equilibrium in the presence of the efficiency recovery mechanism  \cite{brun2013convergence}.
In principle, each strategic firm take turns to update its bid by solving its individual bidding problem \eqref{eq:strategic_bidding_modified} based on the current bids from other firms, and the process iterates until no firm makes changes to its bid. The detailed steps include:

\noindent
\textbf{Step 1}: All firms initialize random positive bids $(b_i(0),i\in\mathcal{N})$, and the clearing price $p$ is set according to \eqref{eq:clearing_price}.

\noindent
\textbf{Step 2}: At $k^{\mathrm{th}}$ iteration with $k=1,2,\dots$, based on the current bids $B_{-i}(k-1)$ of other firms, firm $i$, with 
\begin{equation}
i = \left\{
    \begin{aligned}
   & |\mathcal{N}| \ , \qquad \qquad  \mathrm{if~}k \mathrm{~~mod~~} |\mathcal{N}|=0 \ , \\
    &    k \mathrm{~~mod~~} |\mathcal{N}| \ ,  \ \ \mathrm{otherwise} \ ,
    \end{aligned}\right.
\end{equation}
chooses its optimal bid 
\begin{equation}
    b_i(k)=\arg\max_{b_i:\eqref{eq:clearing_price}} \ {\pi_{i}^{\dagger}(b_i;B_{-i}(k-1))}  \ .
\end{equation}


\noindent
\textbf{Step 3}: The bids of other firms carry over to this iteration:
\begin{equation}
    b_j(k) = b_j(k-1) \ , \ \ j\neq i \ .
\end{equation}
With all these bids $B(k)$, the market updates the clearing price at 
\begin{equation}
 p(k)=\frac{d}{\sum_{i}b_{i}(k)} \ . 
\end{equation}

\noindent
\textbf{Step 4}: Increase $k$ by 1 and go back to Step 2 unless $k \mathrm{~~mod~~} |\mathcal{N}|=0$ holds. Check whether all the bids $B$ remain unchanged in the last $|\mathcal{N}|$ iterations. If so, the algorithm has converged; otherwise, increase $k$ by 1 and go back to Step 2.


The best response bidding algorithm in general guarantees convergence to a Nash equilibrium of a well-defined game, as also validated empirically by our extensive numerical tests for the supply function bidding game in the presence of the efficiency recovery mechanism. A rigorous proof for the convergence of such an algorithm in our setup is an ongoing topic.

\subsection{Evaluation Metrics}

The impact of the strategic behavior of individual participants on market outcomes is associated with their intrinsic market power and usually reflected through their ability to manipulate clearing prices \cite{david2001market}.
From above, strategic bidding does deteriorate global market efficiency, which can be fully recovered with the proposed mechanism. However, the imposed subsidy or tax alters the surplus allocation among firms, thus still maintaining local impacts on individuals. To capture such impacts, we propose two evaluation metrics.

The first is the net profit of each firm, i.e., the objective function $\pi^\dagger$ in \eqref{eq:strategic_bidding_modified}. 
This metric is particularly meaningful when the market is self-sufficient such that it is directly comparable with the surplus allocation at the competitive equilibrium, given the same aggregate net profit. The re-allocation points to who benefits from the mechanism to recover efficiency.

The second is a modified version of the Lerner index \cite{lerner1995concept}, which dates back to 1934 and is a widely recognized measure of market power for individual firms. The standard Lerner Index is defined to capture the difference between the clearing price and a firm's marginal cost:
\begin{equation}\label{eq:std_lerner}
    L^{\mathrm{std}}_{i} = 1- \frac{C_{i}^{'}(q_i)}{p} \ .
\end{equation}
The larger this index is, the more market power a firm possesses.
We adopt a similar idea but make minor changes to adapt to our setting, due to the presence of subsidy or tax. In particular, we integrate the extra subsidy or tax into the cost of each individual firm and compute its marginal cost accordingly. On this basis, we propose our modified Lerner index as
\begin{eqnarray}\label{eq:mod_lerner}
       L_{i} = 1- \frac{C_{i}^{'}(q_i) - \frac{ \mathrm{d}\Delta \pi^\dagger}{\mathrm{d} q_i}}{p} \ .
\end{eqnarray}
It in general has the same indication function as its standard version.


\begin{corollary}\label{cor:eq_Lerner}
In the presence of the efficiency recovery mechanism, the modified Lerner index of each strategic firm $i$ at the Nash equilibrium is explicitly given by
\begin{equation}\label{eq:eq_Lerner}
    L_i = \frac{q_i^{*}}{d - q_i^*} \ , \ \ i\in\mathcal{N} \ , 
\end{equation}
where $(q_i^*,i\in\mathcal{N})$ is an optimal solution to \eqref{eq:CE_opt}.
\end{corollary}
The corollary exploits Theorem~\ref{teo:main_NE} and provides an insight into the underlying relation between the equilibrium Lerner index and the efficient dispatch of supply.
Indeed, the more cost effective a firm is in production, the higher its index \eqref{eq:eq_Lerner} is, indicating that the efficiency recovery mechanism does favor such firms.



\section{Case Studies}

In this section, we validate our equilibrium analysis and evaluate the proposed mechanism with the two metrics, individual net profit and the modified Lerner index, in terms of varying $\phi$ and heterogeneity of firms.
The notion of heterogeneity is captured by the diverse cost functions of firms. We adopt quadratic cost functions, parameterized by quadratic and linear coefficients only.
We will mainly use a 3-firm case as an illustrative example. The standard cost functions of the three firms are shown in Table~\ref{Table.Price and system cost}, but we may vary some of them later to study the impact of heterogeneity. The default value for $\phi$ is set to that in \eqref{eq:self-sufficient_phi}.
The inelastic demand is fixed at $d=10,000$ MW.
\begin{table}[thb]
\centering
\caption{Quadratic cost functions of firms}
\begin{tabular}{ccc}
\toprule 
Firm & \tabincell{c}{quadratic coefficient\\($\$/\mathrm{MW}^2$)} & \tabincell{c}{linear coefficient\\($\$/\mathrm{MW}$)}\\
\midrule 
1 &  0.001 & 10\\
2 &  0.005 & 10\\
3 &  0.005 & 10\\
\bottomrule 
\label{Table.Price and system cost}
\end{tabular}
\end{table}

\subsection{Equilibrium Validation}

\begin{figure}[thb]
\centering
\includegraphics[scale=0.65]{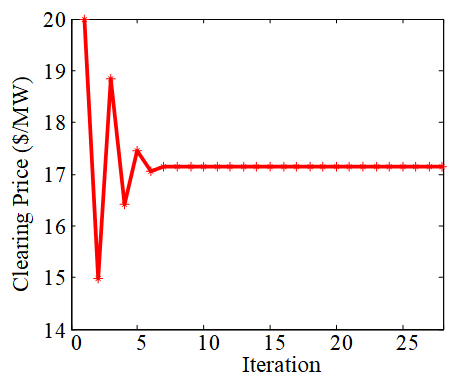}
\caption{Convergence of market clearing price under best response bidding algorithm}
\label{Fig.Equilibrium Clearing price under our proposed mechanism}
\end{figure}

\begin{figure}[thb]
\centering
\includegraphics[scale=0.7]{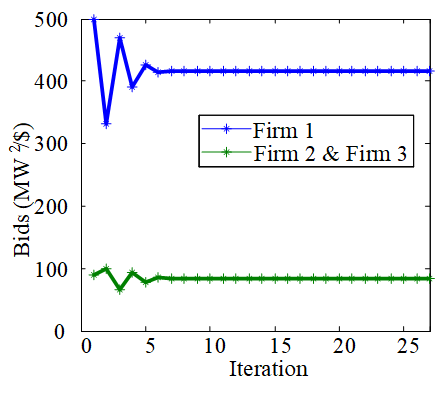}
\caption{Convergence of individual bids under best response bidding algorithm}
\label{Fig.Equilibrium bids under our proposed mechanism}
\end{figure}

\begin{figure}[thb]
\centering
\includegraphics[scale=0.6]{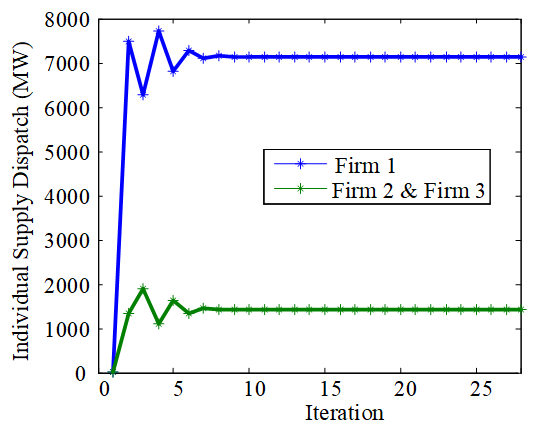}
\caption{Convergence of individual supply dispatch under best response bidding algorithm}
\label{Fig.Equilibrium bidded quantity under our proposed mechanism}
\end{figure}

We run the best response bidding algorithm on the standard 3-firm case to iteratively achieve an equilibrium in the presence of the efficiency recovery mechanism. The convergence of the algorithm is demonstrated in Figure~\ref{Fig.Equilibrium Clearing price under our proposed mechanism}, Figure~\ref{Fig.Equilibrium bids under our proposed mechanism}, and Figure~\ref{Fig.Equilibrium bidded quantity under our proposed mechanism}, showing respectively the evolution of the market clearing price, individual bids, and individual supply dispatch.   


In Table~\ref{Table.Three equilibrium}, we list the three equilibria with the same firms' specifications, where the Nash equilibrium and the competitive equilibrium are directly obtained by solving \eqref{eq:NE_opt} and \eqref{eq:CE_opt}, respectively.
It can be observed that the equilibrium we arrive at using the best response bidding algorithm is consistent with the competitive equilibrium, which is indeed efficient and conforms with Theorem~\ref{teo:main_NE}.
Compared with the Nash equilibrium, the most cost effective firm 1 does not withhold its bid to exaggerate its cost and thus the clearing price is driven back to reflect the system marginal cost in the presence of the efficiency recovery mechanism. 

\begin{table}[thb]
\centering
\caption{Comparison among Nash equilibrium (NE), competitive equilibrium (CE), efficient Nash equilibrium (ENE)}
\begin{tabular}{cccc}
\toprule 
\quad & NE & CE & ENE \\
\midrule 
Social cost ($\$$) & 152627.6 & 135714.3 & 135714.3 
\\
Price ($\$$/MW) & 43.3 & 17.1 & 17.1 \\
Bid 1 ($\mathrm{MW}^2/\$$) & 93.1 & 416.7 
& 416.7 \\
Bid 2 ($\mathrm{MW}^2/\$$) & 68.8 & 83.3 & 83.3 \\
Bid 3 ($\mathrm{MW}^2/\$$) & 68.8 & 83.3 & 83.3 \\
\bottomrule 
\label{Table.Three equilibrium}
\end{tabular}
\end{table}

\begin{table}[thb]
\centering
\caption{Modified Lerner index at ENE}
\begin{tabular}{cccc}
\toprule 
\quad & Firm 1 & Firm 2 & Firm 3 \\
\midrule 
Index  & 2.5 & 0.17 & 0.17 \\
\bottomrule 
\label{Table. Lerner index table}
\end{tabular}
\end{table}

\begin{figure}[thb]
\centering
\includegraphics[scale=0.7]{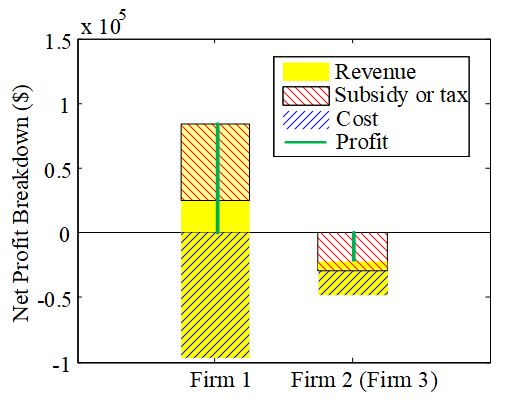}
\caption{Net profit break down at ENE}
\label{Fig.net profit and power market}
\end{figure} 


Table~\ref{Table. Lerner index table} shows the modified Lerner index for each firm at the efficient Nash equilibrium. As Corollary~\ref{cor:eq_Lerner} indicates, the index reflects the truthful cost efficiency of firms, with firm 1 dominating the other two firms, which is consistent with their cost specifications in Table~\ref{Table.Price and system cost}.
Figure~\ref{Fig.net profit and power market} depicts the individual net profit breakdown of each firm. 
We can observe that firms 2 and 3 are taxed to subsidize firm 1, so as to drive all of them to bid truthfully. However, this self-sufficient case suggests that the internal effort to recover efficiency may hurt individual profitability, which discourages market participation in the long run. To avoid such situations, the market may have to conservatively set a small $\phi$, i.e., using external subsidy. This is an intrinsic tradeoff between self-sufficiency and individual profitability for this efficiency recovery mechanism.

\subsection{Impact of $\phi$}

In this subsection, we assess the impact of the constant $\phi$, controlled by the market, in the efficiency recovery mechanism on market surplus allocation among firms. 
Note that Corollary~\ref{cor:eq_Lerner} suggests that the modified Lerner index is independent of $\phi$. Therefore, we will mainly use individual net profit as the metric.
\begin{figure}[thb]
\centering
\includegraphics[scale=0.55]{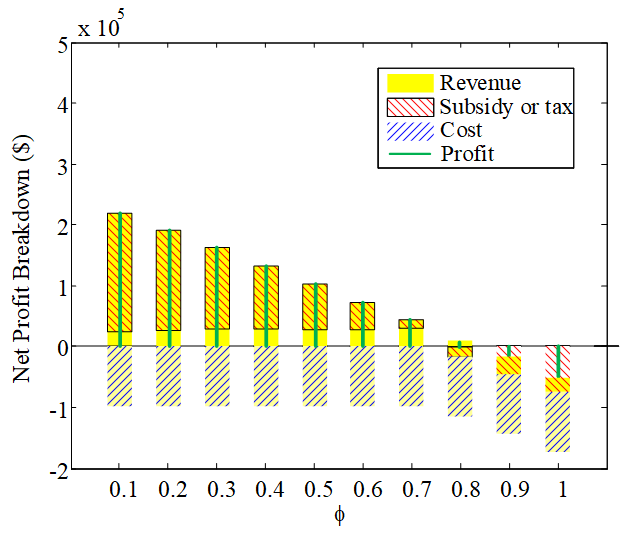}
\caption{Net profit breakdown of firm 1 with varying $\phi$}
\label{Fig.firm 1 net profit with phi}
\end{figure} 
\begin{figure}[thb]
\centering
\includegraphics[scale=0.65]{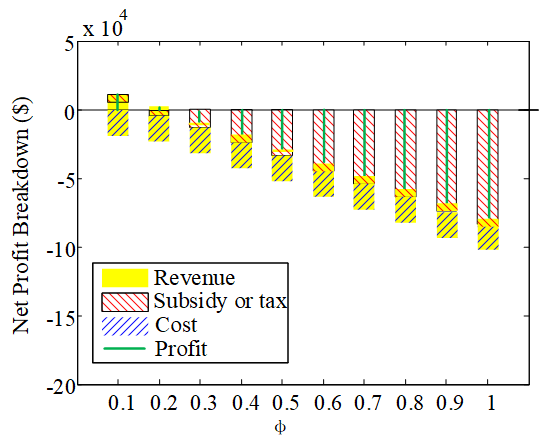}
\caption{Net profit breakdown of firm 2 \& firm 3 with varying $\phi$}
\label{Fig.firm 2/3 net profit with phi}
\end{figure} 

In particular, we vary $\phi$ from 0 to 1 and display the explicit breakdown of the three firms' net profit in Figures~\ref{Fig.firm 1 net profit with phi} and \ref{Fig.firm 2/3 net profit with phi}, including revenue from the market, individual production cost, and subsidy or tax. Note that the net profit of a firm equals its revenue minus cost plus subsidy (or minus tax).
Obviously $\phi$ controls subsidy/tax - in general a firm incurs tax with $\phi$ approaching 1 and enjoys subsidy with $\phi=0$. However, due to the dominance of firm 1 from its cost efficiency, it only requires approximately $\phi\le 0.7$ to enjoy subsidy. In the meantime, the other two firms are much more likely to incur tax, even with only $\phi\ge 0.2$.  
Their different transitioning $\phi$'s between subsidy and tax are not necessarily consistent with the self-sufficient $\phi$ in \eqref{eq:self-sufficient_phi}. Therefore, it remains an open problem to set $\phi$ properly to balance conflicting objectives.





\subsection{Impact of Firm Heterogeneity}

\begin{figure}[thb]
\centering
\includegraphics[scale=0.5]{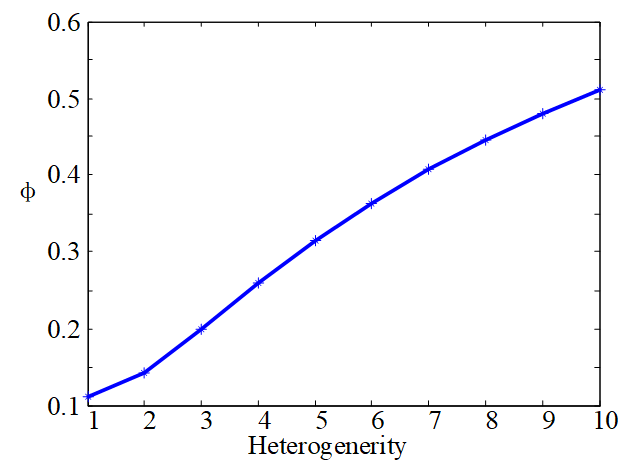}
\caption{Self-sufficiency $\phi$}
\label{Fig.phi}
\end{figure}

In this subsection, we assess the impact of heterogeneity among firms on equilibrium market outcomes in the presence of the efficiency recovery mechanism. In order to capture various degrees of heterogeneity, we vary the quadratic coefficient of firm 1 from 0.005 to 0.0005, and define this coefficient ratio between firm 2/3 and firm 1 as an index for heterogeneity. Therefore, this index of 1 implies a homogeneous case: the more it deviates from 1, the more heterogeneous the case becomes.
Note that as we change the cost specifications, the self-sufficient $\phi$ also varies to guarantee overall subsidy-tax balance, as shown in Figure~\ref{Fig.phi}.

\begin{figure}[thb]
\centering
\includegraphics[scale=0.6]{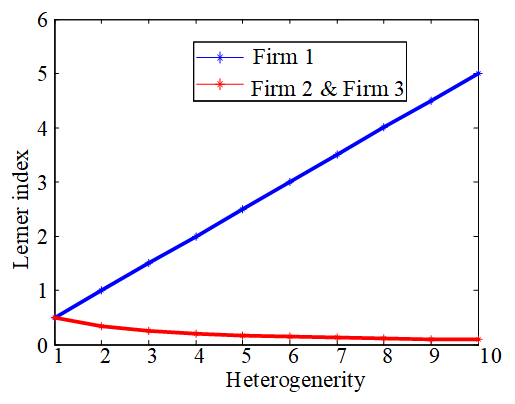}
\caption{Modified Lerner index with varying heterogeneity}
\label{Fig.Market power heterogeneous}
\end{figure}

We first show in Figure~\ref{Fig.Market power heterogeneous} the variation of the modified Lerner index with respect to the heterogeneity index. Starting from the homogeneity case, the Lerner index diverges between firm 1 and firm 2/3. As firm 1 becomes increasingly cost efficient and dominant, the Lerner index of firm 2/3 decreases and converges to zero, while that of firm 1 increases in an approximately linear fashion. This linearity is indeed consistent with Corollary~\ref{cor:eq_Lerner} in the current case of quadratic costs.

\begin{figure}[thb]
\centering
\includegraphics[scale=0.55]{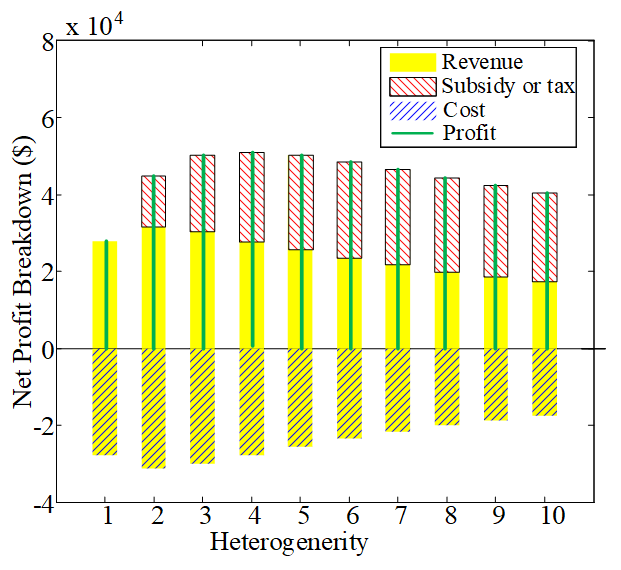}
\caption{Net profit breakdown of firm 1 with varying heterogeneity}
\label{Fig.Profit firm1 heterogeneous}
\end{figure} 

\begin{figure}[thb]
\centering
\includegraphics[scale=0.5]{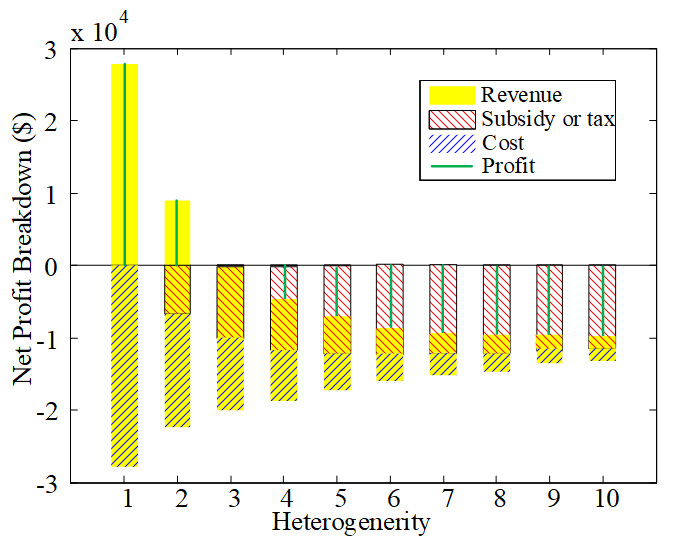}
\caption{Net profit breakdown of firm 2 \& firm 3 with varying heterogeneity}
\label{Fig.profit firm2 heterogeneous}
\end{figure} 

We further show in Figures~\ref{Fig.Profit firm1 heterogeneous} and \ref{Fig.profit firm2 heterogeneous} the net profit breakdown of firm 1 and firm 2/3, respectively, with varying heterogeneity. The increasing heterogeneity implies the dominance of firm 1, for which it has to be incentivized to take on more supply. However, this does not necessarily always earn it more net profit, since its reduced cost coefficient may bring down its production cost and also improves the overall market efficiency with a lower system marginal cost.
Remarkably, we can observe some turning points of heterogeneity where firm 1 obtains the most profit, revenue, cost, and subsidy, respectively.
In the meantime, firms 2-3 takes on less supply that reduces their production cost monotonically. They start, from zero subsidy/tax in the homogeneous case, to incur tax as heterogeneity grows. We can observe similar turning-point behavior for these two firms in their profit, revenue, and subsidy.




\section{Conclusion}

We have studied the competition among firms to meet given inelastic demand in an oligopolistic market with supply function bidding.
To avoid the efficiency loss in social cost caused by firms' strategic bidding behavior, we design an efficiency recovery mechanism that exploits only the market quantities to enforce extra subsidy or tax on individual firms.
The mechanism provably drives the resulting Nash equilibrium to be efficient, and is able to guarantee self-sufficiency with proper tuning of its parameter.
We further evaluate numerically the impact of this mechanism on market outcomes with varying tuning parameters and firm heterogeneity, based on two metrics of individual net profit and a modified version of Lerner index.
Case studies show that the tuning parameter controls the subsidy-tax balance and its choice involves an intrinsic tradeoff between self-sufficiency and individual profitability of the market. 
Besides, more cost efficient firms that enhance their dominance in the market do not necessarily earn more profit as they meanwhile contribute to improving overall efficiency. 

Future extensions are broad. To name a few, first, the market clearing model can be more practical to account for physical constraints, e.g., production capacity, transmission congestion. Second, demand may be elastic and sensitive to prices. Third, it remains a challenge to set the tuning parameter $\phi$ of the mechanism in real operation. We would aim to address them in follow-up studies.




\bibliographystyle{ieeetr}
\bibliography{main}

\end{document}